\documentclass[12pt, a4paper]{article}
\usepackage[T1]{fontenc}
\usepackage[french]{babel}
\usepackage{amsmath,amsthm}
\usepackage{amssymb}
\usepackage{amscd}
\usepackage{euscript}
\usepackage[latin1]{inputenc}
\usepackage{amsfonts}
\usepackage[matrix,arrow,curve]{xy}

\usepackage{cyrillic}

\newcommand{\cyrrm}{\fontencoding{OT2}\selectfont\textcyrup}
\newcommand{\cyrit}{\fontencoding{OT2}\selectfont\textcyrit}

\theoremstyle{plain}
\newtheorem{theo}{Th\'eor\`eme}%[section]
\newtheorem{lem}[theo]{Lemme}
\newtheorem{prop}[theo]{Proposition}
\newtheorem{cor}[theo]{Corollaire}

\theoremstyle{definition}

\newtheorem{theosd}[theo]{Th\'eor\`eme}
\newtheorem{rem}[theo]{Remarque}
\newtheorem{propsd}[theo]{Proposition}

\newcommand{\FF}{\ensuremath{\mathbb{F}}}
\newcommand{\nH}{\ensuremath{H^3_{\mathrm{nr}}}}
\newcommand{\QZl}{\ensuremath{\mathbb Q_l/\mathbb Z_l}}
\newcommand{\ZZl}{\ensuremath{\mathbb Z}/l}

\newcommand{\HH}{\ensuremath{\mathbb H}}
\newcommand{\KK}{\ensuremath{\mathcal K}}
\newcommand{\QZ}{\ensuremath{\mathbb Q/\mathbb Z}}
\newcommand{\Hn}{\ensuremath{H_{\mathrm{nr}}}}
\newcommand{\RH}{\ensuremath{\mathrm{R\underline{Hom}}}}

\newcommand{\eqdef}{\ensuremath{\stackrel{\mathrm{def}}{=}}}

\title{Sur la cohomologie non ramifiée en degré trois d'un produit}
\author{Alena Pirutka}

\begin{document}
\maketitle
\date{}

\begin{abstract}
Soit $\FF$ un corps fini et soit $C$ une courbe projective et lisse sur $\FF$. On établit la nullité du troisième groupe de cohomologie non ramifiée du produit $X\times C$ pour certaines surfaces projectives et lisses $X$ sur $\FF$. Cela s'applique en particulier aux surfaces géométriquement rationnelles.
\end{abstract}

\section{Introduction}
Soit $\FF$ un corps fini. Soit $V/\FF$ une variété projective et lisse, géométriquement intègre, de dimension $3$,  munie d'un morphisme $f:V\to C$ vers une courbe projective et lisse sur $\FF$. Dans un article récent \cite{CTK}, Colliot-Thélène et Kahn ont établi un lien entre la conjecture de Colliot-Thélène et Sansuc sur les zéro-cycles de degré $1$ sur la fibre générique $V_{\eta}$ de $f$ (cf. \cite{CT99}) et le groupe de cohomologie non ramifiée $\nH(V,\mathbb Q_l/\mathbb Z_l(2))$. Plus précisément :
 \begin{quote}
 \textit{Supposons que la conjecture de Tate vaut pour les diviseurs sur $V$ et que le groupe $\nH(V, \QZl(2))$ est divisible. S'il existe sur la surface $V_{\eta}$ une famille de zéro-cycles locaux de degré $1$, orthogonale au groupe de Brauer de $V_{\eta}$ via l'accouplement de Brauer-Manin (cf. \cite {Man71}), alors  il existe sur $V_{\eta}$ un zéro-cycle de degré premier à $l$.}
  \end{quote}

  Par conséquent, on s'intéresse à savoir si le groupe $\nH(V, \QZl(2))$ est nul  pour une telle variété $V$. Plus généralement, on conjecture (\cite{CTK} Conjecture 5.6), que  le troisième groupe de cohomologie non ramifiée s'annule pour une $\FF$-variété projective et lisse, de dimension $3$, géométriquement uniréglée. Pour une variété fibrée en coniques au-dessus d'une surface, c'est un théorème de Parimala et Suresh \cite{PS}. Pour $V$ comme ci-dessus, le premier cas  à examiner est celui d'une fibration triviale $V=X\times C$ où $X$ est une surface géométriquement rationnelle sur $\FF$. Dans cette note on établit  la conjecture pour de telles variétés $V$.
\vspace{0.3cm}

  Pour $k$ un corps on note $\bar k$ une clôture séparable de $k$. Si $X$ est une $k$-variété, on écrit $\bar X=X\times_k\bar k$. Pour $L$ un corps contenant $k$, on écrit
 $X_L= X\times_kL$.
L'énoncé principal est le suivant.
\clearpage

\theosd\label{mXxC}{\textit{Soit $\FF$ un corps fini de caractéristique $p$. Soient $C$ et $X$ des variétés projectives et lisses,  géométriquement intègres, définies sur $\FF$, de dimensions respectives $1$ et $2$.  Soit  $K$ le corps des fonctions de $C$.  Faisons  les hypothèses  :
\begin{itemize}
\item[$\mathrm{(H1)}$] $H^1(\bar X, \mathcal O_{\bar X})=0$;
\item[$\mathrm{(H2)}$] $b_2(\bar X)-\rho(\bar X)=0$;
\item[$\mathrm{(H3)}$] $NS(\bar X)$ est sans torsion;
\item[$\mathrm{(H4)}$]  $A_0(X_{\bar K})=0$.
\end{itemize}
Alors $$\nH(X\times C, \QZl(2))=0$$ pour tout nombre premier $l\neq p$. }\\}

Les hypothèses du théorème sont vérifiées si $\bar X$ est une surface rationnelle. Elles sont aussi satisfaites  pour les surfaces $K3$ supersingulières au sens de Shioda (cf. \cite{F}, ces surfaces n'existent qu'en caractéristique positive). % ainsi que pour les surfaces de Barlow (\cite{Ba}), qui sont des surfaces de type général.

\vspace{0.5cm}

Pour établir le théorème \ref{mXxC}, on procède par diverses réductions. On montre d'abord que le groupe de Chow des $0$-cycles de degré zéro sur $X_K$ est nul, au moins à la $p$-torsion près. On en déduit que tout élément $\xi$ du groupe $$\nH(X\times C, \mu_{l^r}^{\otimes 2})\subset H^3(\FF(X\times C), \mu_{l^r}^{\otimes 2})$$ provient de $H^3_{\acute{e}t}(X_K,  \mu_{l^r}^{\otimes 2})$. Ceci est fait dans la section $\ref{zc}$. Dans la section \ref{cr}, on vérifie que pour les cas que l'on considère, les applications résidus sont compatibles à des applications bord dans la suite spectrale de  Leray. On se ramène ainsi à considérer le groupe  $H^3_{\acute{e}t}(X\times C, \ZZl)$, ce qui nous permet de déduire le résultat dans la section \ref{fp}.

\vspace{1cm}

\subsection{Rappels et notations}

Si $A$ est un groupe abélien et $n$ est un entier, on note $A[n]$ le sous-groupe de A
formé par les éléments annulés par $n$. Pour $l$ un nombre premier on note $A\{l\}$ le sous-groupe de A
formé par les éléments de torsion $l$-primaire.

\'Etant donn\'es un corps $k$ et un entier $n$ inversible sur $k$, on note $\mu_{n}$ le $k$-sch\'ema en groupes (\'etale) des racines $n$-i\`emes de l'unit\'e. Pour $j$ un entier positif, on note $\mu_{n}^{\otimes j}=\mu_{n}\otimes\ldots\otimes\mu_{n}$ ($j$ fois). Lorsque $k$ contient une racine primitive $n$-ième de l'unité, on a un isomorphisme   $\mu_{n}^{\otimes j}\stackrel{\sim}{\to}\mathbb{Z}/n$ pour tout $j$.

Pour $X$ un schéma on note $\mathbb G_m$ le groupe multiplicatif sur $X$ et le faisceau \'etale ainsi d\'efini. On \'ecrit $\mathrm{Br}\,X=H^2_{\acute{e}t}(X,\mathbb G_m)$ pour le groupe de Brauer cohomologique de $X$.

Pour $X$ une $k$-variété propre, on note $A_0(X)$ le groupe de Chow des $0$-cycles de degré zéro sur $X$. On note $\rho(\bar X)=\mathrm{rk}(\mathrm{NS}\,(\bar X))$ le rang du groupe de Néron-Severi de $\bar X$ (cf. \cite{SGA6} XIII 5.1) et $b_i(\bar X)=\mathrm{dim}_{\mathbb Q_l}H^i(\bar X,\mathbb Q_l)$, $l\neq\mathrm{car}\, k$, les nombres de Betti.
\paragraph{Cohomologie non ramifiée.}
Pour $k$ un corps,
 $F$ un corps de fonctions sur $k$,  $n$ un entier  inversible sur $k$,   $i\geq 1$ un entier naturel  et $j\in \mathbb Z$ un entier relatif on définit
\begin{equation*}
\Hn^i(F/k, \mu_n^{\otimes j})=\bigcap\limits_A\mathrm{Ker}[H^i(F, \mu_n^{\otimes j})\stackrel{\partial_{A}}{\to}H^{i-1}(k_A, \mu_n^{\otimes j-1})].
\end{equation*}
\vspace{-0.4cm}

Dans cette formule, $A$ parcourt les anneaux de valuation discr\`ete de rang un, de  corps des fractions $F$, contenant le corps $k$.  Le corps r\'esiduel d'un tel anneau $A$ est not\'e $k_A$ et l'application $\partial_{A}$ est l'application r\'esidu.

Pour $X$ une $k$-vari\'et\'e int\`egre, on note
$\Hn^i(X, \mu_n^{\otimes j})\eqdef\Hn^i(k(X)/k, \mu_n^{\otimes j}),
$
où $k(X)$ est le corps des fonctions de $X$.
On utilise aussi les groupes  $\Hn^{i}(X, \QZ(j))$ (resp.
 $\Hn^{i}(X,\mathbb Q_l/\mathbb Z_l(j))$ pour $l$ un nombre premier) obtenus par passage à la limite inductive.

\paragraph{Rappels de $K$-th\'eorie.}
Pour $X$ un sch\'ema noeth\'erien et $j$ un entier positif on note $\KK_j$ le faisceau de Zariski associ\'e au pr\'efaisceau $U\mapsto K_j(H^0(U, \mathcal O_U))$, le groupe $K_j(A)$ \'etant celui associ\'e par Quillen  \`a l'anneau $A$.
Lorsque $X$ est une vari\'et\'e lisse sur un corps $k$, la conjecture de Gersten, \'etablie par Quillen, permet de calculer les groupes de cohomologie de Zariski $H^i(X, \KK_j)$ comme les groupes de cohomologie du complexe de Gersten. Pour $j=2$ ce complexe s'\'ecrit
\vspace{-0.2cm}
\begin{equation}\label{cxK}K_2k(X)\stackrel{d_2}{\to}\bigoplus\limits_{x\in X^{(1)}}k(x)^*\stackrel{d_1}{\to} \bigoplus\limits_{x\in X^{(2)}} \mathbb Z,\end{equation}
\vspace{-0.4cm}

\noindent où l'application $d_2$ est donn\'ee par le symbole mod\'er\'e et l'application $d_1$ est obtenue par la somme des fl\`eches diviseurs apr\`es normalisation des vari\'et\'es consid\'er\'ees; le groupe $K_2k(X)$ co\"incide avec le groupe  de $K$-th\'eorie de Milnor $K_2^Mk(X)$, quotient de $k(X)^{*}\otimes_{\mathbb Z} k(X)^*$ par le sous-groupe engendr\'e par les \'el\'ements $a\otimes b$ avec $a+b=1$. On en déduit que pour $X$ une vari\'et\'e lisse sur un corps $k$ on a une fl\`eche naturelle
$\mathrm{Pic}(X)\otimes  k^*\to H^1(X, \KK_2).$

\paragraph{Corps globaux.} Si $K$ est un corps global, on note $\Omega$ l'ensemble des places de $K$. On note $K_v$ le complété de $K$ en une place $v$.
Pour $S$ un $K$-tore, on dispose d'une application de restriction $H^i(K,S)\to \prod_{v\in \Omega} H^i(K_v, S)$. On note $$\mbox{\cyrrm{Sh}}^i(K,S)=\mathrm{Ker}[H^i(K,S)\to \prod_{v\in \Omega} H^i(K_v, S)].$$

\subsection{Conséquences des hypothèses}
Les propriétés suivantes sont utilisées plusieurs fois dans la suite.
\prop\label{consh}{Soit $\bar X$ une variété projective et lisse, définie sur un corps séparablement clos $\bar k$. Soit $n>0$ un entier, $(n, \mathrm{car}\,\bar k)=1$. Supposons que $\bar X$  vérifie les hypothèses  $\mathrm{(H1)}- \mathrm{(H3)}$ du théorème \ref{mXxC}. Alors
 \begin{itemize}
 \item[(i)] le groupe $\mathrm{Pic}\,\bar X$ est de type fini sans torsion;
 \item[(ii)] si $\bar k\subset\bar K$ une extension de corps séparablement clos, la flèche naturelle $\mathrm{Pic}\,\bar X\to \mathrm{Pic}\,\bar X_{\bar K}$ est un isomorphisme;
 \item[(iii)] $H^1_{\acute{e}t}(\bar X,\mu_n)=0$;  on a de plus
 $H^3_{\acute{e}t}(\bar X,\mu_n)=0$ si $\bar X$ est une surface;
 \item[(iv)] $Br\,\bar {X}\{l\}=0$ pour tout $l\neq \mathrm{car}\,\bar k$;
  \item[(v)]  on a un isomorphisme $\mathrm{Pic}\,\bar X/n\stackrel{\simeq}{\to} H^2_{\acute{e}t}(\bar X,\mu_n)$;
 \item[(vi)] le noyau et le conoyau de l'application naturelle $\mathrm{Pic}\,\bar X\otimes \bar k^*\to H^1(\bar X, \KK_2)$, sont des groupes uniquement divisibles par tout entier premier à $\mathrm{car}.\,\bar k$.
 \end{itemize}}
 \proof{ L'hypothèse $\mathrm{(H1)}$ implique qu'on a un isomorphisme $\mathrm{Pic}\,\bar X\stackrel{\sim}{\to} NS(\bar X)$ %La dimension de  l'expáce vectoriel $H^1(\bar X, \mathcal O_{\bar X})$ est égale à la dimension de la variété de Picard de $X$.
  où le groupe  $NS(\bar X)$ est un groupe de type fini, et il est sans torsion d'après $\mathrm{(H3)}$, d'où l'énoncé $(i)$. Ce groupe  est invariant par changement des corps séparablement clos car c'est le groupe des composantes du schéma de Picard de $\bar X$, d'où $(ii)$. La suite exacte de Kummer $$0\to\mu_n\to \mathbb G_m\to\mathbb G_m$$ permet d'identifier le groupe $H^1_{\acute{e}t}(\bar X,\mu_n)$ avec le groupe $\mathrm{Pic}\,\bar X[n]$ qui est nul d'après $(i)$. Si $X$ est une surface, cela implique que  $H^3_{\acute{e}t}(\bar X,\mu_n)=0$ par dualité de Poincaré (\cite{Mi} VI.11.2). L'énoncé $(iv)$ est une conséquence des hypothèses $\mathrm{(H2)}$ et $\mathrm{(H3)}$ d'après \cite{GB} II 3.1 et III. 8.3. On déduit alors $(v)$ de la suite de Kummer. L'énoncé $(vi)$  est plus délicat, il est démontré par Colliot-Thélène et Raskind \cite{CTR}, qui utilisent aussi les conjectures de Weil et les résultats de Merkurjev et Suslin \cite{MS}, \cite{Su}.  Soient $K$ le noyau et $C$ le conoyau de l'application $\mathrm{Pic}\,\bar X\otimes \bar k^*\to H^1(\bar X, \KK_2)$.  Sous l'hypothèse $\mathrm{(H2)}$, les énoncés  2.12 et 2.14 de \cite{CTR} impliquent que  $K$ est un groupe uniquement divisible par tout entier premier à $p$ et que $C$ est une somme directe d'un groupe uniquement divisible par tout entier premier à $p$ et du groupe $\oplus_l H^3(\bar X,\mathbb Z_l(1))\{l\}$,  groupe qui est nul d'après l'hypothèse $\mathrm{(H3)}$ (\cite{GB} II 3.1 et III 8.3).  Notons que si $X$ est une surface géométriquement rationnelle, on peut établir $(vi)$ plus facilement (cf. \cite{B} 1.4).   \qed\\}

 \vspace{0.2cm}

\section{Zéro-cycles sur $X_K$ et réduction à $H^3_{\acute{e}t}(X_K, \mu_{l^r}^{\otimes 2})$}\label{zc}
\prop\label{pa0}{Soit $\FF$ un corps fini de caractéristique $p$. Soient $C$ et $X$ des variétés projectives et lisses,  géométriquement intègres, définies sur $\FF$, de dimensions respectives $1$ et $2$. Soit  $K$ le corps des fonctions de $C$.  Supposons que la surface $\bar X/\bar \FF$ vérifie les hypothèses $\mathrm{(H1)-(H4)}$ du théorème \ref{mXxC}.
  Alors le groupe de Chow $A_0(X_K)$ des zéro-cycles de degré zéro sur $X_K$ est nul, au moins à la $p$-torsion près.  \\}

\proof{Soit $G=Gal(\bar \FF/\FF)$ le groupe de Galois absolu de $\FF$.  D'après la proposition \ref{consh}, le $G$-module $\mathrm{Pic}\,\bar X$ est  de type fini sans torsion. Soit $S$ le $\FF$-tore dual.

D'après les estimations de Lang-Weil \cite{LW}, la surface géométriquement intègre $X$ possède un zéro-cycle de degré $1$. On dispose alors d'un torseur \textit{universel} $\mathcal T$ sur $X$ sous $S$ (\cite{CTSan87}, 2.2.2 et 2.3.4). \`A un tel torseur on associe  un homomorphisme
$$A_0(X_k)\stackrel{\Phi_k^{\tau}}{\to} H^1(k,S_k)$$
défini pour toute extension  (non nécessairement finie) $k/\FF$. Pour $x$ un point fermé de $X_k$ de corps résiduel $k(x)$ on pose
$\Phi_k^{\tau}(x)=Cor_{k(x)/k}([\mathcal T_{k(x)}]),$ où $[\mathcal T_{k(x)}]$ est la classe du torseur $\mathcal T\times_X k(x)$ dans $H^1(k(x), S_{k(x)})$. Le fait que cette application passe au quotient par l'équivalence rationnelle est démontré dans \cite{CTSan77} Prop. 12, p. 198.
Cette construction est fonctorielle pour les extensions des corps, d'oùun diagramme commutatif
\begin{align}
\begin{aligned}
\label{d}\xymatrix{
 A_0(X_K)\ar^{\Phi^{\tau}_K}[d]\ar[r]\ar^{\Psi}[dr]& \prod_vA_0(X_{K_v})\ar^{\prod\Phi^{\tau}_{K_v}}[d]&\\
 H^1(K,S_K)\ar[r] & \prod_vH^1(K_v,S_{K_v}).
}
\end{aligned}
\end{align}

L'image de l'application verticale de droite du diagramme est nulle. En effet, soit $L$ une extension finie de $K_v$ et soit $i_{x_v} : \mathrm{Spec}\,L \to X\times_{\FF}K_v$ l'inclusion du point fermé $x_v$ de $X_{K_v}$ de corps résiduel  $L$. Soit $\mathcal O_{L}$ l'anneau des entiers de $L$ et soit $\kappa$ le corps résiduel de $\mathcal O_{L}$. Le morphisme $i_{x_v}$ se prolonge en un morphisme $$\mathrm{Spec}\,\mathcal O_{L} \to X\times_{\FF}\mathrm{Spec}\, \mathcal O_{L}.$$  D'après la construction,  $\Phi^{\tau}_{K_v}(x_v)$ provient alors de la classe du torseur $\mathcal T\times_{\FF} \mathrm{Spec}\, \mathcal O_{L}$ dans $H^1_{\acute{e}t}(\mathcal O_{L},S\times_{\FF}\mathrm{Spec}\,{\mathcal O_{L}})=H^1(\kappa, S_{\kappa})$, le groupe qui est nul car $\kappa$ est un corps fini. Ainsi
\begin{equation} \label{nul}\Psi(A_0(X_K))=0.\end{equation}

Comme la surface $X_K$ vérifie encore les hypothèses $\mathrm{(H1)-(H4)}$,  l'application $\Phi^{\tau}_K$ coïncide avec l'application $\Phi_B: A_0(X_K)\to H^1(X_K, S_{X_K})$ définie par Bloch \cite{B}, au moins à la $p$-torsion près. Pour les surfaces géométriquement rationnelles sur un corps parfait,  c'est le Théorème 3 de \cite{CTSan81}. En utilisant les résultats de \cite{CTR}, on montrera   dans l'appendice de cette note que ce théorème reste vrai, au moins à la $p$-torsion près, sous les conditions $\mathrm{(H1)-(H4)}$.

Pour la suite de la preuve, on raisonne à la $p$-torsion près. Soit ${\mathfrak{G}}=Gal(\bar K/K)$.
On a une suite exacte (\cite{B}, \cite{CTR} 3.6)
\begin{equation*}H^1({\mathfrak{G}}, K_2\bar K(X_{\bar K})/H^0(X_{\bar K},\mathcal K_2))\to A_0(X_K)\stackrel{\Phi_B}{\to} H^1(K,S_{X_K}).
\end{equation*}
Sous l'hypothèse que $X_K$ possède un zéro-cycle de degré $1$, hypothèse qui est satisfaite ici, le théorème  Hilbert 90 pour $\mathcal K_2$ implique que le groupe $H^1({\mathfrak{G}}, K_2\bar K(X_{\bar K}))$ est nul (\cite{CT83}, corollaire 1 et remarque 5.2 pour le cas de la caractéristique positive). De plus, sous les hypothèses $\mathrm{(H1)}$ et $\mathrm{(H3)}$ de la proposition, le groupe $H^0(X_{\bar K},\mathcal K_2)$ est uniquement divisible d'après \cite{CTR} 1.8. On en déduit que le groupe $ H^1({\mathfrak{G}}, K_2\bar K(X_{\bar K})/H^0(X_{\bar K},\mathcal K_2))$ est nul. L'application $\Phi_B=\Phi^{\tau}_K$ est donc injective.

Dans le diagramme (\ref{d}), le noyau $\mbox{\cyrrm{Sh}}^1(K,S_K)$ de l'application horizontale du bas est nul d'après le lemme ci-dessous. Ainsi l'application $\Psi$ est injective. En comparant avec (\ref{nul}), on obtient le résultat.\qed\\}

\lem{Soit $\FF$ un corps fini, soit $C/\FF$ une courbe projective et lisse, géométriquement intègre et soit $K$ son corps des fonctions. Soit $S$ un  $\FF$-tore. Alors $\mbox{\cyrrm{Sh}}^1(K,S_K)=0.$}
\proof{ D'après \cite{CTSan77} p.199, on dispose d'une résolution flasque de $S$ :
\begin{equation}\label{st} 0\to F \to P\to S\to 0 \end{equation}
où $F$ est un tore flasque et $P$ est quasi-trivial (i.e. son module des caractères est un module de permutation). Sur un corps dont le groupe de Galois est cyclique, tout module flasque est un facteur direct d'un module de permutation (\cite{EM}, \cite{CTSan77}). Le tore $F$ est donc un facteur direct d'un tore quasi-trivial. La suite
(\ref{st}) donne alors une inclusion $\mbox{\cyrrm{Sh}}^1(K,S_K)\hookrightarrow \mbox{\cyrrm{Sh}}^2(K,F_K)$ car $H^1(K, F_K)=0$ d'après le théorème Hilbert $90$. Il suffit donc de montrer que $\mbox{\cyrrm{Sh}}^2(K,R_{L/K}\mathbb G_m)=0$ pour $L/K$ une extension finie. En utilisant le lemme de Shapiro, il suffit de montrer que   $\mbox{\cyrrm{Sh}}^2(L,\mathbb G_m)=0$, ce qui résulte du principe de Hasse pour les algèbres centrales simples sur $L$ (cf. par exemple \cite{GS} 6.5.4). \qed}

\rem{L'énoncé de la proposition n'est plus vrai si $K$ est un corps de fonctions de deux variables sur $\FF$, même si $X/\FF$ est une surface géométriquement rationnelle. En effet,  soit $P\in X(\FF)$, soit $\eta$ le point générique de $X$ et soit $\eta'\in X_{\FF(X)}(\FF(X))$ le point rationnel correspondant. On dispose d'un accouplement $Br\,X_{\FF(X)}\times A_0(X_{\FF(X)})\to Br\,\FF(X)$, qui coïncide avec l'inclusion $Br\,X\hookrightarrow Br\, \FF(X)$ sur $Br\,X\times([\eta']-[P])$. Si le groupe $Br\,X$ n'est pas trivial, on en déduit que $[\eta']-[P]$ est un élément non nul de $A_0(X_{\FF(X)})$. \\}

\cor\label{red1}{Soit $\FF$ un corps fini de caractéristique $p$. Soient $C$ et $X$ des variétés projectives et lisses,  géométriquement intègres, définies sur $\FF$, de dimensions respectives $1$ et $2$. Soit $K$ le corps des fonctions de $C$. Supposons que la surface $\bar X/\bar \FF$ vérifie les hypothèses $\mathrm{(H1)-(H4)}$ du théorème \ref{mXxC}.  Soit $l$ un nombre premier, $l\neq p$ et soit $r\geq 1$ un entier. L'application naturelle $$H^3_{\acute{e}t}(X_K,  \mu_{l^r}^{\otimes 2})\to \nH(\FF(X\times C)/K, \mu_{l^r}^{\otimes 2})$$ est surjective.}\proof{De la suite spectrale de Bloch-Ogus  (cf. \cite{BO})
$$ E_2^{pq}= H^p(X_K, \mathcal H^q_{X_K}(\mu_{l^r}^{\otimes 2}))\Rightarrow H^{p+q}_{\acute{e}t}(X_K,\mu_{l^r}^{\otimes 2})$$
on déduit la suite exacte\begin{equation}\label{sebo}
H^3_{\acute{e}t}(X_K,\mu_{l^r}^{\otimes 2})\to \nH(\FF(X\times C)/K, \mu_{l^r}^{\otimes 2})\to CH^2(X_K)/l^r\stackrel{c}{\to} H^4_{\acute{e}t}(X_K,\mu_{l^r}^{\otimes 2}).\end{equation}

La surface géométriquement intègre $X$ possède un zéro-cycle de degré $1$ (cf. \cite{LW}), on a donc de même pour $X_K$. Comme $A_0(X_K)=0$ d'après la proposition précédente, on déduit que le morphisme de degré $CH^2(X_K)\stackrel{\mathrm{deg}}{\to}\mathbb Z$ est un isomorphisme.
Ainsi l'application composée $$CH^2(X_K)/l^r\to H^4_{\acute{e}t}(X_K,\mu_{l^r}^{\otimes 2})\to H^4_{\acute{e}t}(X_{\bar K},\mu_{l^r}^{\otimes 2}),$$ qui s'identifie au morphisme $\mathrm{deg}/l^r$, est injective. On en déduit que l'application $c$ du diagramme (\ref{sebo}) est injective, d'où le resultat.

\qed\\}

\section{\large{Compatibilité des résidus et réduction à $H^3_{\acute{e}t}(X\times C, \mu_{l^r}^{\otimes 2})$}}\label{cr}
\subsection{Notations et énoncés}\label{ne}
Soit $\FF$ un corps fini de caractéristique $p$.  Soient $C$ et $X$ des variétés projectives et lisses,  géométriquement intègres, définies sur $\FF$, de dimensions respectives $1$ et $2$. Soit $K=\FF(C)$. Supposons que la surface $\bar X/\bar \FF$ vérifie les hypothèses $\mathrm{(H1)-(H4)}$ du théorème \ref{mXxC}.
D'après la section précédente, tout élément de $\nH(\FF(X\times C)/K, \mu_{l^r}^{\otimes 2})$ provient d'un élément de $H^3_{\acute{e}t}(X_K,  \mu_{l^r}^{\otimes 2})$. Le but de cette section est d'établir que tout élément du sous-groupe $\nH(\FF(X\times C)/\FF, \mu_{l^r}^{\otimes 2})$ de $\nH(\FF(X\times C)/K, \mu_{l^r}^{\otimes 2})$ provient d'un élément de $H^3_{\acute{e}t}(X\times C, \mu_{l^r}^{\otimes 2})$.

Avec les notations ci-dessus, on a une suite spectrale
\begin{equation}\label{ssl}
E_2^{pq}=H^p(K, H^q_{\acute{e}t}(X_{\bar K}, \mu_{l^r}^{\otimes 2}))\Rightarrow  H^{p+q}_{\acute{e}t}(X_K, \mu_{l^r}^{\otimes 2}).
\end{equation}
D'après la proposition \ref{consh}, on a $H^1_{\acute{e}t}(X_{\bar K},\mu_{l^r}^{\otimes 2})=H^3_{\acute{e}t}(X_{\bar K},\mu_{l^r}^{\otimes 2})=0$ et on a des isomorphismes $\mathrm{Pic}\,\bar X/l^r(1)\stackrel{\simeq}{\rightarrow }\mathrm{Pic}\,X_{\bar K}/l^r(1)\stackrel{\simeq}{\rightarrow }H^2_{\acute{e}t}(X_{\bar K},\mu_{l^r}^{\otimes 2})$.
On dispose alors d'une flèche de bord  \begin{equation}\label{fbord}d^{1,2}: H^3_{\acute{e}t}(X_K, \mu_{l^r}^{\otimes 2})\to H^1(K, \mathrm{Pic}\,\bar X/l^r(1)),\end{equation} qui est un isomorphisme.

Soit $A$ un anneau de valuation discr\`ete de rang un, de  corps des fractions $K$, contenant le corps $\FF$. On a un modèle $X\times_{\FF}\mathrm{Spec}\,A$ de $X_K$ sur $A$ dont on note $X_{k_A}=X\times_{\FF} k_A$ la fibre spéciale. Dans cette section on établit le résultat suivant.\\

\propsd\label{compres}{\textit{Avec les notations ci-dessus,  il existe une application $$d: H^2(X_{k_A}, \mu_{l^r})\to H^0(k_A, \mathrm{Pic\,\bar X}/l^r )$$
telle que le diagramme suivant
\begin{align}
\begin{aligned}
\label{eqres}\xymatrix{ H^3_{\acute{e}t}(X_K, \mu_{l^r}^{\otimes 2})\ar[r]\ar[d]^{d^{1,2}} & H^2_{\acute{e}t}(X_{k_A}, \mu_{l^r})\ar[d]^d&\\
H^1(K, \mathrm{Pic}\,\bar X/l^r(1))\ar[r] &  H^0(k_A, \mathrm{Pic\,\bar X}/l^r )
}
\end{aligned}
\end{align}
où les flèches horizontales sont des applications résidus, est commutatif. \\}}
\rem{\begin{itemize}\item[(i)] On rappelle la construction des résidus dans la section \ref{ar}.
\item[(ii)] Plus précisément, l'application $d$ de la proposition est la composée de l'application naturelle $$H^2_{\acute{e}t}(X_{k_A}, \mu_{l^r}) \to H^0(k_A, H^2_{\acute{e}t}( X_{\bar k_A},\mu_{l^r}))$$
et de l'identification
$H^2_{\acute{e}t}(X_{\bar k_A},\mu_{l^r}) = \mathrm{Pic}\,(\bar X)/l^r$.\\
\end{itemize}}

\subsection{Une application de bord dans la suite spectrale de  Leray}\label{sHS}

Dans ce paragraphe, on décrit la flèche (\ref{fbord}) à l'aide des techniques de \cite{S}.

Soit $S$ un schéma et soit $V$ un $S$-schéma propre. Soit $\Lambda=\mathbb Z/l^r$ où $l$ est inversible sur $V$ et soit $D^+(V, \Lambda)$ (resp. $D^+(S, \Lambda)$) la catégorie dérivée des complexes bornés inférieurement des faisceaux étales de $\Lambda$-modules sur $V$ (resp. sur $S$). Soit $p: V\to S$ le morphisme structural. Soit $\mathcal F$ un faisceau de $\Lambda$-modules  sur $V$ et soit $F=Rp_*\mathcal F$. On a une suite spectrale de Leray $E_2^{ab}=H^a_{\acute{e}t}(S, R^bp_*\mathcal F)\Rightarrow  H^{a+b}_{\acute{e}t}(V, \mathcal F)$.  Supposons  \begin{equation}\label{h03}H^0_{\acute{e}t}(S, R^3p_*\mathcal F)=0.\end{equation} Sous cette hypothèse, on dispose d'une flèche de bord \begin{equation}\label{bord}
\partial :H^{3}_{\acute{e}t}(V, \mathcal F)\to H^1_{\acute{e}t}(S, R^2p_*\mathcal F),
\end{equation}  qui s'obtient comme suit.

On dispose d'un triangle exact \begin{equation}\label{t1}\tau_{\leq 2}F\to F\to \tau_{\geq 3}F.\end{equation}
On a
\vspace{-0.3cm}
 \begin{equation}\label{nulh}\HH^2(S,\tau_{\geq 3}F)=0\;\mathrm{ et }\;  \HH^3(S,\tau_{\geq 3}F)=0.\end{equation}
La première égalité résulte du fait que  le complexe $\tau_{\geq 3}F$ est concentré en degrés strictement supérieurs à $2$. La seconde résulte de l'hypothèse $(\ref{h03})$ et de l'identification $\HH^3(S,\tau_{\geq 3}F)=H^0_{\acute{e}t}(S, H^3(F))$ provenant de la suite exacte longue d'hypercohomologie associée au triangle
$\tau_{\leq 3}(\tau_{\geq 3}F)\to \tau_{\geq 3}F \to \tau_{\geq 4}F$ où $\HH^3(S, \tau_{\geq 4}F)=0$ et $\tau_{\leq 3}(\tau_{\geq 3}F)=H^3(F)[-3]$.

On obtient alors un isomorphisme
$$\varphi: \HH^3(S,\tau_{\leq 2}F)\stackrel{\simeq}{\to} \HH^3(S,F)$$ de la suite exacte longue d'hypercohomologie associée au triangle $(\ref{t1})$.

On a une application naturelle $\psi:\HH^3(S,\tau_{\leq 2}F)\to H^1_{\acute{e}t}(S, H^2(F))$ provenant du triangle
$\tau_{\leq 1}F\to \tau_{\leq 2}F\to H^2(F)[-2]$.

On construit ainsi une application \begin{equation}\label{appl}\psi\circ\varphi^{-1} : \HH^3(S,F)\to H^1_{\acute{e}t}(S, H^2(F))\end{equation} où $\HH^3(S,F)=H^{3}_{\acute{e}t}(V, \mathcal F)$ et $H^1_{\acute{e}t}(S, H^2(F))= H^1_{\acute{e}t}(S, R^2p_*\mathcal F)$.
D'après \cite{S} Appendix B, cette application coïncide avec l'application de bord (\ref{bord}).

\subsection{Applications résidu}\label{ar}
 Soit $S$ un schéma. Soit $V$ un $S$-schéma, soit $j:U\hookrightarrow V$ un ouvert de $V$ et soit $i:Z\to V$ le fermé complémentaire.  Soit $\Lambda=\mathbb Z/l^r$ où $l$ est inversible sur $V$ et soit $F\in Ob D^+(V, \Lambda)$. On a une suite exacte de faisceaux étales sur $V$$$0\to j_!\Lambda\to\Lambda\to i_*\Lambda\to 0,$$
qui donne un triangle exact dans  $D^+(V, \Lambda)$
$$\RH(i_*\Lambda,F)\to \mathrm{R\underline{Hom}}(\Lambda,F)\to \RH(j_!\Lambda,F). $$
En utilisant l'isomorphisme d'adjonction \cite{SGA4} XVIII 3.1.9.10, on obtient
\begin{itemize}
\item[]$\RH(i_*\Lambda, F)=Ri_*\RH(\Lambda, Ri^!F)=i_*Ri^!F$, car $Rf_*=Rf_!$ pour $f$  propre et $f_*$ est exact pour $f$ un morphisme fini; %(voir \cite{Mi}  II.3.6)
\item[]
$\RH(\Lambda,F)=F;$
\item[]$\RH(j_!\Lambda,F)=Rj_*\RH(\Lambda, Rj^!F)=Rj_*j^*F$ car $j_!=Rj_!$ et $Rj^!=j^*$  pour une immersion ouverte.
    \end{itemize}
On obtient ainsi un triangle exact
\begin{equation} \label{sloc} i_*Ri^!F\to F\to Rj_*j^*F.
\end{equation}

Soit $F=\mu_{l^r}^{\otimes j}$ où, plus généralement, $F$ un faisceau de $\Lambda$-modules qui provient d'un faisceau localement constant sur $S$.  Supposons que $V$ et $Z$ sont lisses et que $Z$ est purement de codimension $1$ dans $V$. Dans ce cas, on dispose d'un morphisme de Gysin $F(-1)[-2]\to Ri^!F$ (cf. \cite{Fuj} 1.2) qui est un quasi-isomorphisme  par pureté (cf.  \cite{SGA4} XVI.3.7).
La flèche $\mathrm{R}\,\Gamma \circ Rj_*j^*F\to \mathrm{R}\,\Gamma \circ i_*Ri^!F[1]$ provenant du triangle $(\ref{sloc})$ donne alors des applications \textit{résidus}
$$H^i(U, \mu_{l^r}^{\otimes j})\to H^{i-1}(Z, \mu_{l^r}^{\otimes (j-1)}).$$

En particulier, cette construction donne les applications résidus du diagramme (\ref{eqres})  (voir aussi \cite{CT} et \cite{Mi} III.1.2).

\subsection{Preuve de la proposition \ref{compres}}
Avec les notations de la proposition \ref{compres}, soit $S=\mathrm{Spec\,}A$. On note $\eta$ son point générique et $s$ le point spécial.   On a un diagramme commutatif
$$\xymatrix{
X_{k_A}\ar[r]^(.35){i}\ar[d]^{p_s}& X\times_{\FF} S\ar[d]^p& X_K\ar[l]_(.35){j}\ar[d]^{p_{\eta}}&\\
 s\ar[r]_{i_S} &S & \eta \ar[l]^{j_S}.
}
$$
Soit $\Lambda=\mathbb Z/l^r$ et soit $\mathcal F=\mu_{l^r}^{\otimes 2}$. On a un diagramme commutatif
\begin{align}
\begin{aligned}\label{dcb}\xymatrix{
0\ar[r]& p^*j_{S!}\Lambda\ar[r]\ar[d]& \Lambda\ar[r]\ar[d]&p^*i_{S*}\Lambda\ar[r]\ar[d]& 0&\\
0\ar[r]& j_{!}\Lambda\ar[r]& \Lambda\ar[r]&i_{*}\Lambda\ar[r]& 0
}
\end{aligned}
\end{align}
obtenu en utilisant l'adjonction de $p_*$ et $p^*$ et les identifications $p_*j_!\Lambda=j_{S!}p_{\eta*}\Lambda=j_{S!}\Lambda$ et $p_*i_*\Lambda=i_{S*}p_{s*}\Lambda=i_{S*}\Lambda$.

On a alors un diagramme commutatif dans $D^+(S, \Lambda)$ :
\begin{align}
\begin{aligned}\label{dd}
\xymatrix{
\RH(j_{!}\Lambda, \mathcal F)\ar[d]\ar[r]&\RH(i_{*}\Lambda, \mathcal F)[1]\ar[d]&\\
\RH(p^*j_{S!}\Lambda, \mathcal F)\ar@{=}[d]\ar[r]&\RH(p^*i_{S*}\Lambda, \mathcal F)[1]\ar@{=}[d]&\\
\RH(j_{S!}\Lambda, Rp_*\mathcal F)\ar[r]&\RH(i_{S*}\Lambda, Rp_*\mathcal F)[1]&\\
\RH(j_{S!}\Lambda, \tau_{\leq 2}Rp_*\mathcal F)\ar[r]\ar[u]^{\iota_1}\ar[d]& \RH(i_{S*}\Lambda,\tau_{\leq 2}Rp_*\mathcal F)[1]\ar[u]^{\iota_2}\ar[d]&\\
\RH(j_{S!}\Lambda, R^2p_*\mathcal F)[-2]\ar[r]& \RH(i_{S*}\Lambda, R^2p_*\mathcal F)[-1].
}
\end{aligned}
\end{align}
Dans ce diagramme, le premier carré est déduit du diagramme $(\ref{dcb})$, le deuxième carré est obtenu par adjonction et les deux derniers carrés proviennent de la construction
du paragraphe \ref{sHS} par fonctorialité.

D'après la proposition \ref{consh}, on a $R^2p_*\mathcal F=\mathrm{Pic}\, \bar X/l^r(1)$ (voir aussi \cite{Mi} VI 8.9). D'après la construction des résidus dans le paragraphe précédent, on a alors que l'application $$H^3_{\acute{e}t}(X_K, \mu_{l^r}^{\otimes 2})\to H^2_{\acute{e}t}(X_{k_A}, \mu_{l^r})$$ (resp. l'application  $H^1(K, \mathrm{Pic}\, \bar X/l^r(1))\to  H^0(k_A, \mathrm{Pic}\, \bar X/l^r))$ du diagramme (\ref{eqres}) provient de la ligne du haut (resp. du bas) du diagramme $(\ref{dd})$, si l'on  applique le foncteur $\mathrm{R}\Gamma$ et si l'on prend ensuite l'hypercohomologie en degré trois.
  Par le même argument que dans le paragraphe \ref{sHS}, après cette opération les morphismes $\iota_1$ et $\iota_2$ du diagramme  $(\ref{dd})$ induisent des isomorphismes. En effet, la proposition \ref{consh} implique que le faisceau $R^3p_*\mathcal F$ est nul; on dispose alors d'un triangle exact $\tau_{\leq 2}Rp_*\mathcal F\to Rp_*\mathcal F\to \tau_{\geq 4}Rp_*\mathcal F$ où le complexe $\tau_{\geq 4}Rp_*\mathcal F$ est concentré en degrés strictement supérieurs à $3$.

Le diagramme $(\ref{eqres})$ provient alors de la première et la dernière lignes du diagramme $(\ref{dd})$.\qed\\

\subsection{Réduction à $H^3_{\acute{e}t}(X\times C, \mu_{l^r}^{\otimes 2})$}
\prop\label{red2}{Soit $\FF$ un corps fini de caractéristique $p$. Soient $C$ et $X$ des variétés projectives et lisses,  géométriquement intègres, définies sur $\FF$, de dimensions respectives $1$ et $2$.   Soit $l$ un nombre premier, $l\neq p$. Si la surface $\bar X/\bar \FF$ vérifie les hypothèses $\mathrm{(H1)-(H4)}$ du théorème \ref{mXxC},  alors tout élément de $\nH(\FF(X\times C)/\FF, \mu_{l^r}^{\otimes 2})$  provient d'un élément de $H^3_{\acute{e}t}(X\times C, \mu_{l^r}^{\otimes 2})$. }
\proof{Soit  $K=\FF(C)$. Soit $\xi$ un élément de $\nH(\FF(X\times C)/\FF, \mu_{l^r}^{\otimes 2})$. Comme $\xi$ est en particulier non ramifié par rapport à $K$,  il provient d'un élément $\xi'$ de $H^3_{\acute{e}t}(X_K,  \mu_{l^r}^{\otimes 2})$ d'après le corollaire \ref{red1}.

On dispose d'un isomorphisme $d^{1,2}: H^3_{\acute{e}t}(X_K,  \mu_{l^r}^{\otimes 2})\stackrel{\sim}{\to} H^1(K, \mathrm{Pic}\,\bar X/l^r(1))$ (cf. $(\ref{fbord})$).
La proposition \ref{compres} implique que  $\xi''=d^{1,2}(\xi')$ est non ramifié, i.e. $$\xi''\in H^1_{nr}(K/\FF, \mathrm{Pic}\,\bar X/l^r(1)).$$

Par un argument analogue à \cite{CT}  4.2.1, on en déduit que $\xi''$ provient du groupe  $H^1_{\acute{e}t}(C, \mathrm{Pic}\,\bar X/l^r(1))$. En effet, $\xi''$ est dans l'image de $H^1_{\acute{e}t}(\mathrm{Spec}\,\mathcal O_{C,c}, \mathrm{Pic}\,\bar X/l^r(1))$ pour tout point $c$ de $C$, ce que l'on peut voir de la suite $(\ref{sloc})$. On en déduit que $\xi''$ provient de  $H^1_{\acute{e}t}(C, \mathrm{Pic}\,\bar X/l^r(1))$ en utilisant  la suite de Mayer-Vietoris (cf.  \cite{CT} 3.8.2).

Notons $\pi:X\times C\to C$ le morphisme de projection. Soit $F= \mu_{l^r}^{\otimes 2}$, vu comme faisceau étale sur $X\times C$. En utilisant la proposition \ref{consh}, on obtient que les faisceaux  $R^1\pi_*F$ et $R^3\pi_*F$ sont nuls car c'est le cas pour leurs  fibres géométriques. Le faisceau $R^2\pi_*F$ s'identifie à  $\mathrm{Pic}\,\bar X/l^r(1)$ (cf. \cite{Mi} VI.8.9). Le groupe $H^4_{\acute{e}t}(C, \mu_{l^r}^{\otimes 2})$ est nul, car $cd_l\,C\leq 3$ pour $C$ une courbe sur un corps fini. La suite spectrale de Leray $E_2^{pq}=H^p(C, R^q\pi_*F)\Rightarrow  H^{p+q}_{\acute{e}t}(X\times C,  \mu_{l^r}^{\otimes 2})$ donne alors un morphisme de bord surjectif $H^3_{\acute{e}t}(X\times C,  \mu_{l^r}^{\otimes 2})\to H^1_{\acute{e}t}(C, \mathrm{Pic}\,\bar X/l^r(1))$. Par fonctorialité, on a le diagramme commutatif suivant

$$\xymatrix{
H^3_{\acute{e}t}(X_K,  \mu_{l^r}^{\otimes 2})\ar[r]_(.43){d^{1,2}}^(0.43){\simeq}& H^1(K, \mathrm{Pic}\,\bar X/l^r(1))&\\
 H^3_{\acute{e}t}(X\times C,  \mu_{l^r}^{\otimes 2})\ar@{->>}[r]\ar[u] & H^1_{\acute{e}t}(C, \mathrm{Pic}\,\bar X/l^r(1))\ar[u].
}
$$

 Puisque $\xi''=d^{1,2}(\xi')$ provient de $H^1_{\acute{e}t}(C, \mathrm{Pic}\,\bar X/l^r(1))$, on en déduit  que $\xi$ provient  de $H^3_{\acute{e}t}(X\times C,  \mu_{l^r}^{\otimes 2})$.  \qed\\}

\section{Fin de la preuve du théorème \ref{mXxC}}\label{fp}

Pour établir le théorème \ref{mXxC}, il suffit  de montrer que le groupe $\nH(X\times C, \mu_l^{\otimes 2})$ est nul. En effet, ce groupe s'identifie au  sous-groupe de $l$-torsion du groupe  $\nH(X\times C, \QZl(2))$, ce qui est une conséquence du théorème de Merkurjev et Suslin (cf. \cite{MS} p. 339).   Puisque le degré de l'extension de $\FF$ obtenue en ajoutant les racines $l$-ièmes de l'unité est premier à $l$, on peut de plus supposer que $\FF$ contient les racines $l$-ièmes de l'unité par un argument de corestriction. D'après la proposition \ref{red2},  tout élément du groupe $\nH(X\times C, \ZZl)$ provient d'un élément de  $H^3_{\acute{e}t}(X\times C, \ZZl)$. Le théorème \ref{mXxC} résulte donc de la proposition suivante. \\

\prop\label{prnul}{Soit $\FF$ un corps fini de caractéristique $p$.  Soient $C$ et $X$ des variétés projectives et lisses,  géométriquement intègres, définies sur $\FF$, de dimensions respectives $1$ et $2$.  Supposons que la surface $\bar X/\bar \FF$ vérifie les hypothèses $\mathrm{(H1)}-\mathrm{(H4)}$  du théorème \ref{mXxC}. Soit $l\neq p$ un nombre premier. Supposons que $\FF$ contient une racine primitive $l$-ième de l'unité.  Alors l'image de $H^3_{\acute{e}t}(X\times C, \ZZl )$ dans $H^3(\FF(X\times C), \ZZl)$ est nulle.}
\proof{Soit $G=Gal(\bar \FF/\FF).$
Pour établir la proposition, on va utiliser deux lemmes suivants.\\}

\lem\label{lsuites}{Supposons les hypothèses de la proposition \ref{prnul} vérifiées. \begin{itemize}\item[(i)] On a alors une suite exacte  où toutes les flèches sont des flèches évidentes : \begin{multline}
 0\to H^3_{\acute{e}t}(X,\QZl(2))\to H^3_{\acute{e}t}(X_{\FF(C)}, \QZl(2))\to\\\to H^3_{\acute{e}t}(X_{\bar \FF(C)}, \QZl(2))^G\to 0.
\end{multline}
\item[(ii)] L'application $\mathrm{Pic}\,\bar X\otimes H^1_{\acute{e}t}(\bar C,\ZZl)\to H^3_{\acute{e}t}(\bar X\times \bar C, \ZZl)$ induite par les cup-produits  est un isomorphisme.
    \end{itemize}}

\proof{Soient $K=\FF(C)\subset \bar \FF(C)$ et $\mathfrak{G}=Gal(\bar K/\bar \FF(C))$ le groupe de Galois absolu du corps $\bar \FF(C)$. La proposition \ref{consh}, appliquée à $\bar X$ et à $X_{\bar K}$, implique que  $H^1_{\acute{e}t}(\bar X,\QZl(2))=H^3_{\acute{e}t}(\bar X,\QZl(2))=0$ et que $H^1_{\acute{e}t}(X_{\bar K},\QZl(2))=0$. On a également $H^3_{\acute{e}t}(\bar C,\QZl(2))=0$ car $cd_l\, \bar C\leq 2$.

 En appliquant la suite spectrale de Leray pour $X_{\FF(C)}\to \mathrm{Spec}\, \FF(C)$, on obtient alors une suite exacte
 \small
 $$0\to H^1(G, H^2_{\acute{e}t}(X_{\bar \FF(C)}, \QZl(2)))\to  H^3_{\acute{e}t}(X_{\FF(C)}, \QZl(2))\to H^3_{\acute{e}t}(X_{\bar \FF(C)}, \QZl(2))^G\to 0.\\$$
 \normalsize
 De même, les suites spectrales de Leray pour $X\to \mathrm{Spec}\, \FF$ et pour $X_{\bar \FF(C)}\to \mathrm{Spec}\, \bar \FF(C)$ donnent
\begin{align*}
H^1(G,H^2(\bar X,\QZl(2)))&\stackrel{\simeq}{\to}  H^3_{\acute{e}t}(X,\QZl(2)),\\
H^2_{\acute{e}t}(X_{\bar \FF(C)}, \QZl(2))&\stackrel{\simeq}{\to}H^2_{\acute{e}t}(X_{\bar K}, \QZl(2))^{\mathfrak{G}}.
\end{align*}
\normalsize
 Pour établir l'énoncé $(i)$ du lemme, il suffit donc de montrer qu'on a un isomorphisme $ H^2_{\acute{e}t}(\bar X, \QZl(2))\stackrel{\simeq}{\to} H^2_{\acute{e}t}(X_{\bar K}, \QZl(2))^{\mathfrak{G}}$. Or ces deux $G$-modules s'identifient à $\mathrm{Pic}\,\bar X\otimes \QZl(1)$. En effet, d'après la propostion \ref{consh} appliquée à $\bar X$ et à $X_{\bar K}$, on a  $\mathrm{Pic}\,\bar X/l^r(1) \stackrel{\simeq}{\to}H^2_{\acute{e}t}(\bar X, \mu_{l^r}^{\otimes 2})$ et $\mathrm{Pic}\, X_{\bar K}/l^r(1)\stackrel{\simeq}{\to}  H^2_{\acute{e}t}(X_{\bar K},\mu_{l^r}^{\otimes 2})$ pour tout $r>0$,
 et la flèche naturelle $\mathrm{Pic}\,\bar X\to \mathrm{Pic}\, X_{\bar K}$ est un isomorphisme.

Montrons $(ii)$. Comme les groupes de cohomologie $H^i_{\acute{e}t}(\bar X,\ZZl)$ et $H^j_{\acute{e}t}( \bar C,\ZZl)$  sont des $\ZZl$-espaces vectoriels, les groupes $\mathrm{Tor}_r^{\ZZl}(H^i_{\acute{e}t}(\bar X,\ZZl), H^j_{\acute{e}t}(\bar C,\ZZl))$ sont nuls pour $r>0$. La formule de Künneth  fournit alors une décomposition \begin{equation*}\bigoplus\limits_{i+j=m} H^i_{\acute{e}t}(\bar X,\ZZl)\otimes_{\ZZl} H^j_{\acute{e}t}(\bar C,\ZZl)\stackrel{\simeq}{\to} H^m_{\acute{e}t}(\bar X\times \bar C, \ZZl)
\end{equation*} où les flèches sont induites par les cup-produits
(cf. \cite{SGA4} XVII.5.4.3).

Pour  $m=3$ on obtient
$H^2_{\acute{e}t}(\bar X,\ZZl)\otimes_{\ZZl}H^1_{\acute{e}t}(\bar C, \ZZl)\stackrel{\simeq}{\to} H^3_{\acute{e}t}(\bar X\times \bar C, \ZZl)$, puisque les groupes $H^1_{\acute{e}t}(\bar X,\ZZl)$, $H^3_{\acute{e}t}(\bar X,\ZZl)$ et $H^3_{\acute{e}t}(\bar C,\ZZl)$ sont nuls.

D'après la proposition \ref{consh},  $\mathrm{Pic}\,\bar X/l\stackrel{\simeq}{\to} H^2_{\acute{e}t}(\bar X,\ZZl)$, d'où
 \begin{equation*} H^3_{\acute{e}t}(\bar X\times \bar C, \ZZl)=\mathrm{Pic}\,\bar X/l\otimes_{\mathbb Z/l}H^1_{\acute{e}t}(\bar C,\ZZl)=\mathrm{Pic}\,\bar X\otimes H^1_{\acute{e}t}(\bar C,\ZZl),
 \end{equation*}
d'où le deuxième énoncé du lemme.
 \qed\\}

 Soit $\iota$ l'application naturelle
$$\iota: [\mathrm{Pic}\,\bar X\otimes H^1_{\acute{e}t}(\bar C,\ZZl)]^G\to [\mathrm{Pic}\,\bar X\otimes H^1_{\acute{e}t}(\bar C,\QZl(1))]^G\to [\mathrm{Pic}\,\bar X\otimes \bar \FF(C)^*\otimes \QZl]^G.$$

Le lemme suivant montre qu'on a une inclusion \begin{equation}\label{incl}\iota[(\mathrm{Pic}\,\bar X\otimes H^1_{\acute{e}t}(\bar C,\ZZl))^G]\subset [\mathrm{Pic}\,\bar X\otimes \bar \FF(C)^*]^G\otimes\QZl.
                                              \end{equation}

\lem\label{perm}{Soit $\FF$ un corps fini, soit $G=Gal(\bar \FF/\FF)$ et soit $l$ un nombre premier, $l\neq car.\FF$.  Soit $C$ une courbe projective et lisse, géométriquement intègre, définie sur $\FF$.
\begin{itemize}
\item[(i)] Pour $P$ un $G$-module de permutation, le groupe  $H^1(G, P\otimes  \bar \FF(C)^*)$ est nul  et
 l'application naturelle $$[P\otimes \bar \FF(C)^*]^G\otimes\QZl\to [P\otimes  \bar\FF(C)^*\otimes \QZl]^G$$ est surjective.
\item[(ii)] Pour $M$ un $G$-module de type fini sans torsion, on a $$H^1(G, M\otimes  H^1_{\acute{e}t}(\bar C, \QZl(1)))=0.$$
\item[(iii)] Pour $M$ un $G$-module de type fini sans torsion, l'image de l'application naturelle  $$ [M\otimes H^1_{\acute{e}t}(\bar C, \QZl(1))]^G\to [M\otimes \bar \FF(C)^*\otimes \QZl]^G$$ est contenue dans le sous-groupe $[M\otimes \bar \FF(C)^*]^G\otimes\QZl$ de $[M\otimes \bar \FF(C)^*\otimes \QZl]^G$.
\end{itemize}
}

\proof{
Pour montrer l'énoncé $(i)$, il suffit de considérer le cas où $P=\mathbb Z[G/H]$ où $H=Gal(\bar \FF/L)$ est un sous-groupe ouvert de $G$.  Pour tout $G$-module $A$ et pour tout $i\geq 0$, on a $H^i(G,\mathbb Z[G/H]\otimes  A)\stackrel{\simeq}{\to} H^i(H,A)$  (cf. \cite{Se} I.2.5 ou \cite{NSW} p.59). On a donc $H^1(G, \mathbb Z[G/H]\otimes  \bar \FF(C)^*)=H^1(H, \bar\FF(C)^*)$. D'après le théorème de Hilbert $90$, on a $H^1(H, \bar \FF(C)^*)= H^1(Gal(\bar \FF(C)/L(C)), \bar\FF(C)^*)=0$.

Le diagramme suivant $$\xymatrix{
 [P\otimes \bar \FF(C)^*]^G\ar[r]\ar[d]^{\times l} &[P\otimes  \bar\FF(C)^*/\bar \FF(C)^{*l^r}]^G\ar[d]^{\times l}&\\
[P\otimes \bar \FF(C)^*]^G\ar[r] &[P\otimes  \bar\FF(C)^*/\bar \FF(C)^{*l^{r+1}}]^G&
}$$
est commutatif. En passant à la limite, pour montrer la deuxième partie de l'énoncé $(i)$, il suffit de montrer que pour tout $r>0$ l'application naturelle $[P\otimes \bar \FF(C)^*]^G\to [P\otimes  \bar\FF(C)^*/\bar \FF(C)^{*l^r}]^G$ est surjective.

 On a une suite exacte $$0\to \mu_{l^r}\to  \bar \FF(C)^*\to \bar \FF(C)^*/\mu_{l^r}\to 0.$$ On déduit que le groupe $H^1(G,P\otimes \bar \FF(C)^*/\mu_{l^r})$ est  nul, car $H^1(G, P\otimes \bar \FF(C)^*)=0$  d'après ce qui précède et $H^2(G, P\otimes \mu_{l^r})=0$ car $cd_l\FF\leq 1$.
La suite exacte de $G$-modules
$$0\to \bar \FF(C)^*/\mu_{l^r}\stackrel{\times l^r}{\to} \bar \FF(C)^*\to \bar \FF(C)^*/\bar\FF(C)^{*l^r}\to 0$$
donne alors la surjection
$[P\otimes \bar \FF(C)^*]^G\to [P\otimes  \bar\FF(C)^*/\bar\FF(C)^{*l^r}]^G.$

Montrons l'énoncé $(ii)$. On a une suite exacte $$0\to I\to P\to M\to 0,$$ où $P$ est un $G$-module de permutation. Cette suite est $\mathbb Z$-scindée et elle reste donc exacte après tensorisation avec le groupe $H^1_{\acute{e}t}(\bar C, \QZl(1))$. Comme $cd_l\,\FF\leq 1$, le groupe $H^2(G, I \otimes  H^1_{\acute{e}t}(\bar C, \QZl(1)))$ est nul,
d'où une application surjective $H^1(G, P\otimes  H^1_{\acute{e}t}(\bar C, \QZl(1)))\to H^1(G, M\otimes  H^1_{\acute{e}t}(\bar C, \QZl(1)))$. Il suffit donc de considérer le cas où $M$ est un module de permutation. Par le même argument que ci-dessus, on se ramène à montrer que $H^1(H, H^1_{\acute{e}t}(\bar C, \QZl(1)) )=0$ pour $H=Gal(\bar\FF/L)$. Cela résulte du théorème de Lang \cite{L} sur les groupes algébriques connexes sur un corps fini. En effet, on a $H^1_{\acute{e}t}(\bar C, \QZl(1))=J(\bar \FF)\{l\}$ où $J$ est la jacobienne de $C$ et $H^1(H,J)=0$ (\textit{loc. cit.}).

 L'énoncé $(ii)$ appliqué à $I$ donne $H^1(G, I\otimes H^1_{\acute{e}t}(\bar C,\QZl(1)))=0$, d'où une application surjective $[P\otimes H^1_{\acute{e}t}(\bar C,\QZl(1))]^G\to [M\otimes H^1_{\acute{e}t}(\bar C,\QZl(1))]^G$. On déduit alors $(iii)$ du diagramme commutatif suivant, dans lequel l'application\\ $[P\otimes \bar \FF(C)^*]^G\otimes\QZl\to [P\otimes  \bar\FF(C)^*\otimes \QZl]^G$  est surjective d'après $(i)$.

\small
 $$\xymatrix{
[P\otimes H^1_{\acute{e}t}(\bar C,\QZl(1))]^G\ar@{->>}[r]\ar[d]& [M\otimes H^1_{\acute{e}t}(\bar C,\QZl(1))]^G\ar[d]& \\
 [P\otimes \bar \FF(C)^*\otimes \QZl]^G\ar[r]& [M\otimes \bar \FF(C)^*\otimes \QZl]^G& \\
[P\otimes \bar \FF(C)^*]^G\otimes\QZl\ar[r]\ar@{->>}[u]& [M\otimes \bar \FF(C)^*]^G\otimes\QZl.\ar[u]&
}
$$

\normalsize
\qed\\}

On passe maintenant à la preuve de la proposition \ref{prnul}.\\

Soient \begin{align*}NH^3_{\acute{e}t}(X_{\FF(C)},\ZZl)&=\mathrm{ker}[H^3_{\acute{e}t}(X_{\FF(C)},\ZZl)\to H^3(\FF(X\times C),\ZZl)],\\
NH^3_{\acute{e}t}(X_{\FF(C)},\QZl(2))&=\mathrm{ker}[H^3_{\acute{e}t}(X_{\FF(C)},\QZl(2))\to H^3(\FF(X\times C),\QZl(2))].\end{align*}

On a  le diagramme commutatif suivant
\footnotesize
\begin{align}
\begin{aligned}\label{d1}\xymatrix{
 H^3_{\acute{e}t}(X\times C, \ZZl)\ar@{->>}[r]\ar[d] &H^3_{\acute{e}t}(\bar X\times\bar C, \ZZl)^G\ar[d]&\\
 H^3_{\acute{e}t}(X_{\FF(C)}, \ZZl)\ar@{->>}[r]\ar[d] &H^3_{\acute{e}t}(X_{\bar \FF(C)}, \ZZl)^G\ar[d]&\\
 H^3_{\acute{e}t}(X_{\FF(C)}, \QZl(2))\ar@{->>}[r] &H^3_{\acute{e}t}(X_{\bar \FF(C)}, \QZl(2))^G&
}
\end{aligned}
\end{align}

\normalsize
\noindent où les applications horizontales sont surjectives car $cd\,G\leq 1$.

Soit $d: H^3_{\acute{e}t}(X\times C, \ZZl)\to H^3_{\acute{e}t}(X_{\bar \FF(C)}, \QZl(2))^G$ l'application diagonale de ce diagramme. Soit $\alpha\in H^3_{\acute{e}t}(X\times C, \ZZl)$. Il suffit de montrer que l'image de $\alpha$ dans $H^3_{\acute{e}t}(X_{\FF(C)}, \ZZl)$ est incluse dans le groupe $NH^3_{\acute{e}t}(X_{\FF(C)}, \ZZl)$.
Le théorème de Merkurjev et Suslin \cite{MS} implique que l'application $H^3(\FF(X\times C), \ZZl)\to H^3(\FF(X\times C), \QZl(2))$ est injective.   Il suffit alors de montrer que l'image de $\alpha$ dans  le groupe $H^3_{\acute{e}t}(X_{\FF(C)}, \QZl(2))$ est incluse dans  $NH^3_{\acute{e}t}(X_{\FF(C)}, \QZl(2))$.

D'après le lemme \ref{lsuites}, le noyau de l'application horizontale du bas du diagramme (\ref{d1}) s'identifie au groupe $H^3_{\acute{e}t}(X, \QZl(2))$. Ce groupe est inclus dans le groupe $NH^3_{\acute{e}t}(X_{\FF(C)}, \QZl(2))$ car  pour $X$   une surface projective et lisse, géométriquement connexe, définie sur un corps fini,  le groupe $H^3_{nr}(X, \QZl(2))$ est nul d'après  \cite{CTSS}, p.790. Pour établir la proposition \ref{prnul}, il suffit alors de montrer que l'on peut relever $d(\alpha)$ en un élément $\delta$ de $NH^3_{\acute{e}t}(X_{\FF(C)}, \QZl(2))$. Pour construire cet élément on procède comme suit.

D'après le lemme \ref{lsuites}, on a un isomorphisme $$(\mathrm{Pic}\,\bar X\otimes H^1_{\acute{e}t}(\bar C,\ZZl))^G\simeq H^3_{\acute{e}t}(\bar X\times\bar C, \ZZl)^G$$ induit par les cup-produits. L'image de $\alpha$ dans $H^3_{\acute{e}t}(\bar X\times\bar C, \ZZl)^G$ provient alors d'un élément $\beta\in  (\mathrm{Pic}\,\bar X\otimes H^1_{\acute{e}t}(\bar C,\ZZl))^G$. D'après (\ref{incl}), on a que $\iota(\beta)$ est un élément de $(\mathrm{Pic}\,\bar X\otimes \bar \FF(C)^*)^G\otimes\QZl$. Soit $\gamma$ son image  dans $H^1(X_{\bar \FF(C)}, \mathcal{K}_2)^G\otimes\QZl$. On a alors le diagramme commutatif suivant :
\small
\begin{align}
\begin{aligned}\label{dl1}
\xymatrix{
 \iota(\beta)\in [\mathrm{Pic}\,\bar X\otimes \bar \FF(C)^*]^G\otimes \QZl\ar@{^{(}->}[r]\ar[d] &[\mathrm{Pic}\,\bar X\otimes \bar \FF(C)^*\otimes \QZl]^G\ar[d]&\\
 \gamma\in H^1(X_{\bar \FF(C)}, \mathcal{K}_2)^G\otimes \QZl\ar[r] &[H^1(X_{\bar \FF(C)}, \mathcal{K}_2)\otimes \QZl]^G.&
}
\end{aligned}
\end{align}
\normalsize

 D'après les résultats de Colliot-Thélène et Raskind \cite{CTR} 4.3, sous les hypothèses de la proposition \ref{prnul} l'application naturelle $H^1(X_{\FF(C)}, \mathcal{K}_2)\to H^1(X_{\bar \FF(C)}, \mathcal{K}_2)^G$ est un isomorphisme à la $car.\FF$-torsion près. On a donc que $\gamma$ provient d'un élément $\gamma'\in H^1(X_{\FF(C)}, \mathcal{K}_2)\otimes\QZl.$ On dispose d'une application $H^1(X_{\FF(C)},\mathcal{K}_2)\otimes \QZl\to NH^3_{\acute{e}t}(X_{\FF(C)}, \QZl(2))$ (voir \cite{CTc} (3.11)). Soit $\delta$ l'image de $\gamma'$ dans ce dernier groupe. On a un diagramme commutatif
\begin{align}
\begin{aligned}\label{dl2}
\xymatrix{
 \gamma'\in H^1(X_{\FF(C)}, \mathcal{K}_2)\otimes \QZl\ar[r]\ar[d] &[H^1(X_{\bar \FF(C)}, \mathcal{K}_2)\otimes \QZl]^G\ar[d]&\\
 \delta\in NH^3_{\acute{e}t}(X_{\FF(C)}, \QZl(2))\ar[r] &H^3_{\acute{e}t}(X_{\bar \FF(C)}, \QZl(2))^G.&
}
\end{aligned}
\end{align}
Soit $x\in \bar X^{(1)}$ et soit $[x]$ sa classe dans $\mathrm{Pic}\,\bar X$. L'image de $[x]\otimes H^1(\bar C,\ZZl)$ par l'application composée $\mathrm{Pic}\,\bar X\otimes H^1_{\acute{e}t}(\bar C,\ZZl)\to H^1(X_{\bar\FF(C)}, \mathcal{K}_2)/l\to H^3(X_{\bar \FF(C)}, \ZZl)$ provient d'un cup-produit avec la classe canonique de $x$ dans $H^2_x(X_{\bar \FF(C)}, \ZZl)$ (cf. \cite{CTc}, \cite{BO}). Ainsi $d(\alpha)$ est l'image de $\beta$ par l'application composée \begin{multline*}[\mathrm{Pic}\,\bar X\otimes H^1_{\acute{e}t}(\bar C,\ZZl)]^G\stackrel{\iota}{\to}[\mathrm{Pic}\,\bar X\otimes \bar \FF(C)^*\otimes \QZl]^G\to\\\to [H^1(X_{\bar \FF(C)}, \mathcal{K}_2)\otimes \QZl]^G\to H^3_{\acute{e}t}(X_{\bar \FF(C)}, \QZl(2))^G. \end{multline*}
D'après les diagrammes $(\ref{dl1})$ et $(\ref{dl2})$, $d(\alpha)$ provient alors de $\delta\in NH^3_{\acute{e}t}(X_{\FF(C)}, \QZl(2))$.
\qed

\subsection*{Appendice : comparaison des applications  $\Phi^{\tau}$ et $\Phi_B$}

Soit $k$ un corps  de caractéristique $p\geq 0$. Soit $X$ une $k$-variété projective et lisse, telle que  la variété $\bar X/\bar k$ vérifie les hypothèses suivantes :
\begin{itemize}
\item[$\mathrm{(H1)}$] $H^1(\bar X, \mathcal O_{\bar X})=0$;
\item[$\mathrm{(H2)}$] $b_2-\rho=0$;
\item[$\mathrm{(H3)}$] $NS(\bar X)$ est sans torsion;
\item[$\mathrm{(H4)}$] $A_0(\bar X)=0$;
\item[$\mathrm{(H5)}$] $X$ possède un zéro-cycle de degré $1$.\\
\end{itemize}

Soit $G=Gal(\bar k/k)$. Soit $S$ le $k$-tore dual du $G$-module  $\mathrm{Pic}\,\bar X$ (cf. proposition \ref{consh}). Sous l'hypothèse $\mathrm{(H5)}$, on dispose d'un torseur universel $\mathcal T$ sur $X$ sous $S$ (\cite{CTSan87}) et on construit  une application
$A_0(X)\stackrel{\Phi^{\tau}}{\to} H^1(k, S)$
(\cite{CTSan77} p.198).

Dans le cas où $X$ est une surface, on dispose également d'une application $\Phi_B: A_0(X){\to} H^1(k, S)$ définie par Bloch \cite{B}. Si de plus $X$ est une surface géométriquement rationnelle définie sur un corps parfait, Colliot-Thélène et Sansuc ont établi que ces deux applications coïncident (\cite{CTSan81}, Thm. 3).
Dans cet appendice, on rappelle leur preuve et on utilise les résultats de Colliot-Thélène et Raskind \cite{CTR} pour nous convaincre que les hypothèses  $\mathrm{(H1)}-\mathrm{(H5)}$ suffisent, au moins à la $p$-torsion près. Par un argument de transfert, cette dernière restriction nous permet de  ne considérer que des cycles $\xi=\sum n_ix_i$ de $A_0(X)$  où  $x_i$ sont des points fermés de $X$ dont les corps résiduels  $k(x_i)$ sont des extensions séparables de $k$. %En générale,  si $L$ est une extension normale finie de $k$ contenant tous les corps $k(x_i)$, alors il existe une extension purement inséparable $k'\subset L$ de $k$ de degré $p^m$ pour certain $m>0$,   telle que $L/k'$ est une extension galoisienne.  Soit $\xi_{k'}$ l'image de $\xi$ dans $X_{k'}$. Si $\Phi_B(\xi_{k'})=\Phi^{\tau}(\xi_{k'})$, alors $\Phi_B(p^m\xi)=\Phi^{\tau}(p^m\xi)$ par un argument de corestriction.

\paragraph*{Application $\Phi_B$.} On suppose que $X$ est une surface.
 On a une application naturelle $\mathrm{Pic}\,\bar X\otimes \bar k^*\to H^1(\bar X, \KK_2)$, dont le noyau et le conoyau sont des groupes uniquement divisibles par tout entier premier à $p$ (\cite{CTR},  cf. proposition \ref{consh} ci-dessus). On en déduit que l'application
\begin{equation}
\beta : H^1(G, \mathrm{Pic}\,\bar X\otimes \bar k^*)\to H^1(G, H^1(\bar X, \KK_2))
\end{equation}
est un isomorphisme à la $p$-torsion près.

Soit $Y\subset X$ un diviseur et soit $U\subset X$ l'ouvert complémentaire. Soit $x\in \bar X^{(1)}$ un point qui n'est pas dans $\bar Y$, soit $D_x$ l'adhérence de $x$ dans $\bar X$ et soit $\pi_x :\tilde D_x\to D_x$ la normalisation de $D_x$. Soit $\bar k(x)_{\bar Y}^*$ le groupe des fonctions rationnelles sur  $D_x$, qui sont inversibles en tout point de $\pi_x^{-1}(\bar Y)$.
On a un complexe (cf. p. \pageref{cxK})
\begin{equation*}K_2\bar k(\bar X)\stackrel{d_2}{\to}\bigoplus\limits_{x\in \bar X^{(1)}}\bar k(x)^*\stackrel{d_1}{\to} \bigoplus\limits_{x\in \bar X^{(2)}} \mathbb Z,\end{equation*}
d'où une suite exacte
\begin{equation}\label{suitek}
 0\to \mathcal Z_{\bar Y}\to \bigoplus\limits_{x\in \bar U^{(1)}}\bar k(x)_{\bar Y}^*\to (\bigoplus\limits_{x\in \bar X}^{0} \mathbb Z)_{\bar Y}\to 0,
\end{equation}
où le terme de gauche (resp. de droite) est défini comme le noyau (resp. l'image) de la restriction de $d_1$  au groupe du milieu.
En particulier, on  a une application $\mathcal Z_{\bar Y}\to H^1(\bar X, \KK_2)$.
En passant à la cohomologie dans la suite (\ref{suitek}), on obtient une application $\partial : H^0(G, (\bigoplus\limits_{x\in \bar X}^{0} \mathbb Z)_{\bar Y})\to H^1(G, \mathcal Z_{\bar Y})$. Par composition avec l'inverse de l'application $\beta$ on obtient alors une application
$\Phi_{B,Y} :  H^0(G, (\bigoplus\limits_{x\in \bar X}^{0} \mathbb Z)_{\bar Y})\to H^1(G, \mathrm{Pic}\,\bar X\otimes \bar k^*)$
définie à la $p$-torsion près.

D'après la construction, le groupe $(\bigoplus\limits_{x\in \bar X}^{0} \mathbb Z)_{\bar Y}$ est inclus dans le groupe des zéro-cycles de degré zéro sur $\bar U$. D'après \cite{CTSan81} p. 426, l'hypothèse $A_0(\bar X)=0$ implique que  l'élément $\xi$  provient d'un zéro-cycle $\xi'$ sur $X$ tel que le cycle $\bar \xi'$ sur $\bar X$ appartient au groupe $(\bigoplus\limits_{x\in \bar X}^{0} \mathbb Z)_{\bar Y}$ construit comme ci-dessus pour certain $Y\subset X$. En effet, on peut écrire $\bar \xi'=\mathrm{div}\,(\sum_{i=1}^r f_i|_{D_i})$ où $D_i\subset \bar X$ sont des courbes intègres et $f_i\in \bar k(D_i)^*$.  Soit $\pi_i: \tilde D_i\to D_i$ la normalisation de $D_i$. On peut alors prendre $U$ un ouvert de $X$ contenant $\pi_i(\mathrm{div}_{\tilde D_i}(f_i))$ et $Y$ son complémentaire. Notons que dans cette construction on peut de plus supposer que $\mathrm{Pic}\,\bar U=0$.  On pose alors
$$\Phi_B(\xi)=\Phi_{B,Y}(\xi').$$

\paragraph*{Application $\Phi^{\tau}$.}   Soit $U\subset X$ un ouvert tel que $\mathrm{Pic}\,\bar U=0$.
On a une suite exacte des $G$-modules
 \begin{equation}\label{gmod}
0\to \bar k[U]^*/\bar k^*\to \mathrm{Div}_{\bar Y}\bar X\to \mathrm{Pic}\,\bar X\to 0.
\end{equation}
Soit
\begin{equation}\label{dgmod}
  1\to S\to M\to R_U\to 1
\end{equation}
la suite duale.

D'après \cite{CTSan87} 2.3.4, la projection $\bar k[U]^*\to \bar k[U]^*/\bar k^*$ admet une section $\sigma$ (en général, cela est assuré par l'existence d'un torseur universel).
La section $\sigma$ donne une application $\phi_{\sigma} : U\to R_U$.
 D'après \cite{CTSan87} 2.3.4, la restriction de $\mathcal T=\mathcal T^{\sigma}$ à $U$ se déduit via $\phi_{\sigma}$ du torseur $M$ défini par la suite $(\ref{dgmod})$.
   Soit $x$  un point fermé de $U$, tel que son corps résiduel $k(x)$ est une extension séparable de $k$.  Alors $\Phi^{\tau}(x)\in H^1(k,S)\stackrel{\sim}{\to} \mathrm{Ext}^1_k(\hat S, \bar k^*)$ est obtenu de la classe de l'extension (\ref{gmod}) via  $$\mathrm{Ext}^1_k(\hat S,\bar k[U]^*/\bar k^*)\stackrel{\sigma}{\to} \mathrm{Ext}^1_k(\hat S,\bar k[U]^*)\stackrel{\mathrm{ev}_x}{\to}\mathrm{Ext}^1_k(\hat S,(k(x)\otimes_k\bar k)^*)\stackrel{N_{k(x)/k}}{\to} \mathrm{Ext}^1_k(\hat S, \bar k^*)$$ (cf. \cite{CTSan81} p. 425).

\paragraph*{Comparaison.} Considérons un accouplement $$\bigoplus\limits_{x\in \bar U^{(1)}}\bar k(x)_{\bar Y}^*\times \mathrm{Div}_{\bar Y}\bar X\to \bar k$$ défini comme suit. Soit $f\in \bar k(x)_{\bar Y}^*$, soit $D$ l'adhérence de $x$ dans $\bar X$ et soit $\pi: \tilde D\to D$ la normalisation.  Soit $\Delta\subset \bar X$  un diviseur irréductible à support dans $\bar Y$ et soit $\pi^{-1}(\Delta)=\sum n_i(y_i)$. On pose $(f, \Delta)=\prod f(y_i)^{n_i}$. D'après \cite{B} A.8, le diagramme
$$\xymatrix{
\bigoplus\limits_{x\in \bar U^{(1)}}\bar k(x)_{\bar Y}^*\times \mathrm{Div}_{\bar Y}\bar X\ar[r]\ar@<-2pc>[d]& \bar k^*\ar@{=}[d]&\\
 (\bigoplus\limits_{x\in \bar X}^{0} \mathbb Z)_{\bar Y} \times \bar k[U]^*/\bar k^*\ar@<-2pc>[u]_{\mathrm{div}}\ar[r]^(.70){ev}& \bar k^*.
}
$$
est commutatif. En effet, soit $g\in \bar k[U]^*/\bar k^*$ et soit $g'$ son image dans $\bar k(x)^*/\bar k^*$. Il s'agit donc de vérifier que $f(\mathrm{div}_D(g'))=g'(\mathrm{div}_D(f))$, ce qui résulte de la loi de réciprocité de Weil.

On obtient ainsi un diagramme commutatif
\begin{footnotesize}
\begin{align}
\begin{aligned}\label{ld1}\xymatrix{
0\ar[r]&  \mathcal Z_{\bar Y}\ar[r]\ar[d] &\bigoplus\limits_{x\in \bar U^{(1)}}\bar k(x)_{\bar Y}^*\ar[r]\ar[d]& (\bigoplus\limits_{x\in \bar X}^{0}\mathbb Z)_{\bar Y}\ar[r]\ar[d]& 0&\\
0\ar[r]& \mathrm{Hom}(\mathrm{Pic}\,\bar X, \bar k^*)\ar[r]& \mathrm{Hom}(\mathrm{Div}_{\bar Y}\bar X, \bar k^*)\ar[r]&\mathrm{Hom}(\bar k[U]^*/\bar k^*, \bar k^*)\ar[r]&0.
}
\end{aligned}
\end{align}
\end{footnotesize}

On a  aussi un triangle commutatif
\begin{align}
\begin{aligned}\label{ld2}\xymatrix{
\mathcal Z_{\bar Y}\ar[r]\ar[dr]& H^1(\bar X, \KK_2)& \mathrm{Pic}\,\bar X\otimes \bar k^*\ar[l]\ar[ld]&\\
& \mathrm{Hom}(\mathrm{Pic}\,\bar X, \bar k^*)&
}
\end{aligned}
\end{align}
où l'application de droite est déduite du produit d'intersection sur $\mathrm{Pic}\,\bar X$
(\cite{B} A.11). En effet, on dispose d'un accouplement $H^1(\bar X, \KK_2)\times \mathrm{Div}_{\bar Y}\bar X\to \bar k^*$ compatible avec le produit d'intersection sur  $\mathrm{Pic}\,\bar X$. Pour $\Delta\in \mathrm{Div}_{\bar Y}\bar X$ une courbe irréductible dont on note  $\tilde \Delta\to X$ la normalisation, l'accouplement est défini par  $H^1(\bar X, \KK_2)\to H^1(\tilde \Delta, \KK_2)\to \bar k^*$ où l'application  $H^1(\tilde \Delta, \KK_2)\to \bar k^*$ est induite par la norme (cf. \textit{loc. cit.}).

On déduit de $(\ref{ld1})$ et $(\ref{ld2})$ le diagramme commutatif suivant (défini à la $p$-torsion près)
$$\xymatrix{
H^0(G,(\bigoplus\limits_{x\in \bar X}^{0}\mathbb Z)_{\bar Y})\ar[r]\ar[d]& H^1(G, \mathcal Z_{\bar Y})\ar[d]^{\gamma}\ar[r]& H^1(G, \mathrm{Pic}\,\bar X\otimes \bar k^*)\ar[dl]^{\simeq}&\\
\mathrm{Hom}_{G}(\bar k[U]^*/\bar k^*, \bar k^*)\ar[r]\ar[r]& H^1(G,\mathrm{Hom}(\mathrm{Pic}\,\bar X, \bar k^*)).
}
$$
D'après la construction, pour tout $\xi\in H^0(G,(\bigoplus\limits_{x\in \bar X}^{0}\mathbb Z)_{\bar Y})\subset A_0(X)$, on a alors que $\Phi_B(\xi)=\Phi^{\tau}(\xi)$.

\paragraph{Remerciements.} La proposition \ref{compres} résulte d'une discussion avec O. Wittenberg; je lui sais gré de m'avoir   expliqué les techniques utilisées pour la preuve. Je voudrais aussi remercier J.-L. Colliot-Thélène  pour de nombreux conseils.

\end{document}